\documentclass[12pt]{amsart}
\usepackage{amsmath,amsfonts,amssymb}

\def\res{\mathop{\rm{res}}}

\def\ssum{\mathop{\sum\!\sum}}

\def\sumb{\sideset{}{^\flat}\sum}

\def\le{\leqslant}

\def\ge{\geqslant}

\renewcommand{\mod}{\mathop{\rm{mod}}}

\def\cA{\mathcal{A}}

\def\cR{\mathcal{R}}

\numberwithin{equation}{section}

\newtheorem{cor}{COROLLARY}[section]
\newtheorem{prop}[cor]{PROPOSITION}
\newtheorem{lem}[cor]{LEMMA}

\theoremstyle{definition}

\newtheorem{theorem}{THEOREM}

\begin{document}
\title{\bf Selberg's sieve of irregular density}

\author[Friedlander]{J.B. Friedlander$^*$}
\thanks{$^*$\ Supported in part by NSERC grant A5123}
\author[Iwaniec]{H. Iwaniec}

\maketitle

\medskip

{\bf Abstract:}
 We study certain aspects of the Selberg sieve, in particular when sifting by rather thin sets of primes. We derive new results for the lower bound sieve suited especially for this setup and we apply them in particular to give a new sieve-propelled proof of Linnik's theorem on the least prime in an arithmetic progression in the case of the presence of exceptional zeros. 

\section{\bf Introduction}

Let $\cA= (a_n)$ be a finite sequence of numbers $a_n\ge 0$. Our goal is to estimate the sifting function
\begin{equation}\label{eq:1.1} 
S(\cA,z) = \sum_{(n,P(z))=1}a_n
\end{equation}
where $P(z)$ is the product of all primes $p<z$. 

To this end we assmue that the congruence sums
\begin{equation}\label{eq:1.2} 
A_d = \sum_{n\equiv 0(\mod d)}a_n
\end{equation}
are well-approximated by
\begin{equation}\label{eq:1.3} 
A_d = g(d)X +r_d \quad {\rm if} \,\, d\mid P(z) 
\end{equation} 
where $g(d)$ is a multiplicative function (called the density) satisfying 
\begin{equation}\label{eq:1.4} 
0\le g(p) <1 
\end{equation} 
and $X>0$ is a fixed number such that the error terms $r_d$ are small 
on average. 

We assume that the ``remainder''
\begin{equation}\label{eq:1.5} 
R(y)= \sum_{d<y}\tau_3(d)|r_d| 
\end{equation} 
is considerably smaller than the contribution to $S(\cA,z)$ from the main 
terms in~\eqref{eq:1.3}, so it is negligible. The larger the value of $y$ 
that can be used, the stronger the estimates that can be obtained. 
Here, the mysterious looking presence of the divisor function $\tau_3$ is 
due to our use of the Selberg sieve weights. 

The $\Lambda^2$-sieve of Selberg yields a lovely upper bound 
\begin{equation}\label{eq:1.6} 
S(\cA,z)\le  XJ^{-1} +R(\Delta^2)
\end{equation} 
where
\begin{equation}\label{eq:1.7} 
J= \sum_{d\le\Delta,\,d|P(z)}h(d) 
\end{equation} 
and $h(d)$ is the multiplicative function with
\begin{equation}\label{eq:1.8} 
h(p)=\frac{g(p)}{1-g(p)}.
\end{equation} 
Obviously, 
\begin{equation}\label{eq:1.9} 
J\le \prod_{p|P(z)}(1+h(p))=V(z)^{-1} 
\end{equation} 
where
\begin{equation}\label{eq:1.10} 
V(z)=\prod_{p|P(z)}(1-g(p)) .
\end{equation} 
Hence, the main term in the upper bound~\eqref{eq:1.6} exceeds $XV(z)$. 
Therefore, the remainder is insignificant if, say,
\begin{equation}\label{eq:1.11} 
R(\Delta^2) \ll X V(z)(\log \Delta)^{-1}. 
\end{equation} 

But, we still need a lower bound for $J$ in~\eqref{eq:1.6}. A strong 
and useful bound is difficult to establish unless we make some further 
assumptions about the density function $g(d)$ or its companion $h(d)$. It is 
quite often the case that $h(d)$ behaves nicely in the sense that the 
Dirichlet series 
\begin{equation}\label{eq:1.12} 
D(s)= \sumb_{d\ge 1}h(d)d^{1-s}
\end{equation} 
admits an analytic continuation to ${\rm Re}\, s >\frac{1}{2}$ , with only 
a pole at $s=1$. 

Here and in what follows, the superscript $\flat$ indicates that the summation 
variable is restricted to squarefree numbers. 

\smallskip 
Chapter 7 of [Opera] provides a plethora of major examples with proofs. In  
Section 7.9 an asymptotic formula for $J$ is given under the assumption
\begin{equation}\label{eq:1.13} 
\sum_{p\le x}g(p)\log p = \kappa \log x +O(1) 
\end{equation} 
for every $x\ge 2$, where $\kappa$ is  a positive constant called the sieve 
dimension. This approximation tells us that $g(p)p$ is $\kappa$ on average 
provided that $p$ is sufficiently large in terms of some defining parameters 
of $g$, so that the error term in~\eqref{eq:1.13} can be ignored. 

However, in some practical sieve applications, the function $g(p)$ appears in 
segments in which it is not uniformly distributed. For example, in Section 24.2 
of [Opera] we have
\begin{equation}\label{eq:1.14} 
\begin{cases}
g(p)p = 1+\chi(p)\bigl(1-\frac{1}{p}\bigr) \quad {\rm if}\,\, p\nmid q \\
g(p) =0 \quad {\rm if}\,\, p\mid q ,
\end{cases} 
\end{equation} 
where $\chi (\mod q)$ is a real, non-principal character. In this 
case~\eqref{eq:1.13} holds with $\kappa =1$ but, so far as we know, 
with a terribly poor error term.   

Fortunately, sieve methods can also 
produce useful estimates when one has access only to upper bounds for $g(p)$.  
In the above example, we can use the trivial bound $g(p)p\ < 2$ giving
\begin{equation}\label{eq:1.15} 
\sum_{p\le x}g(p)\log p \le 2\log x +O(1). 
\end{equation} 
To get an absolutely bounded error term we have here sacrificed the 
sieve dimension so the results are weaker. Still, if $z$ is small in the 
logarithmic scale (as in the Fundamental Lemma of sieve theory) compromising 
$\kappa$ does not significantly affect the output. 

Even if $z$ is relatively large there can be significant consequences in 
the case of upper bounds. Unfortunately, the sieve of superficially enlarged 
dimension may yield a negative lower bound for $S(\cA,z)$ in a range of $z$ 
where it is expected to be positive. 

In these notes we are concerned with $g(p)p$ fluctuating unpredictably 
within the segment 
\begin{equation}\label{eq:1.16} 
0\le g(p)p < 2 .
\end{equation} 
We shall establish a positive lower bound
\begin{equation}\label{eq:1.17} 
S(\cA,z)\gg XV(z)
\end{equation} 
for $z$ quite large, provided that $g(p)p$ is small on average. To this end 
we could go through the recurrence formula of Buchstab via the Fundamental 
Lemma. However, we can derive very explicit and neat results by means of 
Selberg's lower-bound sieve method.

\section{\bf Selberg's Lower-Bound Sieve}

Following Selberg [S], we have 
\begin{equation}\label{eq:2.1} 
S(\cA,z)\ge S^-(\cA,z)
\end{equation} 
where
\begin{equation}\label{eq:2.2} 
S^-(\cA,z)=\sum_na_n\bigl(1-\sum_{p|n,\, p<z}1\bigr)\bigl(\sum_{d|n}\rho_d\bigr)^2
\end{equation} 
with any real numbers $\rho_d$, $\rho_1=1$. We assume that $\rho_d$ are 
supported on squarefree $d\le \Delta$. Opening the square and applying the 
approximations~\eqref{eq:1.3} we obtain
\begin{equation}\label{eq:2.3} 
S^-(\cA,z)=XW+ \cR 
\end{equation} 
where
\begin{equation}\label{eq:2.4} 
W = \sum_{d_1}\sum_{d_2}\rho_{d_1}\rho_{d_2}
\bigl(g([d_1,d_2])-\sum_{p<z}g([p,d_1,d_2])\bigr)
\end{equation} 
and $\cR $ is the corresponding remainder
\begin{equation}\label{eq:2.5} 
\cR = \sum_{d_1}\sum_{d_2}\rho_{d_1}\rho_{d_2}
\bigl(r_{[d_1,d_2]}-\sum_{p<z}r_{[p,d_1,d_2]}\bigr) .
\end{equation} 

In the quadratic form $W$ in the variables $\rho_d$ we make a linear 
change of variables, specifically setting
\begin{equation}\label{eq:2.6} 
y_d = \frac{\mu(d)}{h(d)}\sum_{m\equiv 0 (\mod d)}g(m)\rho_{m}. 
\end{equation} 
Applying M\"obius inversion, we find
\begin{equation}\label{eq:2.7} 
\rho_{\ell} = \frac{\mu(\ell)}{g(\ell)}\sum_{d\equiv 0 (\mod \ell)}g(d)y_d . 
\end{equation} 
The quadratic form~\eqref{eq:2.4} in the new variables becomes
\begin{equation}\label{eq:2.8} 
W= \sum_dh(d)y_d^2-\sum_{p<z}g(p)\sum_{(d,p)=1}h(d)(y_d-y_{pd})^2 ;
\end{equation} 
see the formula (7.109) in [Opera]. 

\smallskip
{\bf Remarks:} We have tacitly assumed that $g(\ell)$, $h(d)$ do not vanish 
to justify the transformations~\eqref{eq:2.6},~\eqref{eq:2.7}. However, 
after obtaining~\eqref{eq:2.8} we no longer need this slight detour.

It is clear that the support conditions for $\rho_d$ and $y_d$ are the 
same, namely
\begin{equation}\label{eq:2.9} 
d\le \Delta, \quad d  \,\, {\rm squarefree}. 
\end{equation} 
We could apply the stronger conditions $d\le \Delta$, $d|P(z)$ which we 
omit for technical simplifications of our arguments. 

\smallskip 

The normalization $\rho_1=1$ becomes
\begin{equation}\label{eq:2.10} 
\sum_dh(d)y_d=1.
\end{equation} 
Now we choose
\begin{equation}\label{eq:2.11} 
y_d=\frac{1}{H}\log\frac{\Delta}{d} 
\end{equation} 
for $d\le \Delta$, $d$ squarefree,
where 
\begin{equation}\label{eq:2.12} 
H=\sumb_{d\le \Delta}h(d)\log\frac{\Delta}{d} .
\end{equation} 
Recall this means we are restricting $d$ to squarefree variables. 

The sieve constituents $\rho_d$ are bounded. Precisely, as in the 
pure $\Lambda^2$-sieve of Selberg, for every squarefree $\ell\le\Delta$ we 
argue as follows:
\begin{equation*}
\begin{aligned}
H & =\sum_{k|\ell}\sum_{\substack{d<\Delta\\(d,\ell)=k}}h(d)\log\frac{\Delta}{d} 
=\sum_{k|\ell}h(k)\sum_{\substack{m<\Delta/k\\(m,\ell)=1}}h(m)\log\frac{\Delta}{mk}\\ 
& \ge\bigl(\sum_{k|\ell}h(k)\bigr) \sum_{\substack{m<\Delta/\ell\\(m,\ell)=1}}h(m)
\log\frac{\Delta}{m\ell}=\frac{h(\ell)}{g(\ell)}
\sum_{\substack{m<\Delta/\ell\\(m,\ell)=1}}h(m)\log\frac{\Delta}{m\ell} .
\end{aligned}
\end{equation*}
On the other hand, pulling the factor $h(\ell)$ out from the 
sum~\eqref{eq:2.7}, we get 
\begin{equation*}
\mu(\ell)\rho_{\ell}= \frac{h(\ell)}{g(\ell)}
\sum_{\substack{m<\Delta/\ell\\(m,\ell)=1}}h(m)\log\frac{\Delta}{m\ell} .
\end{equation*}
Combining these results we find that
\begin{equation}\label{eq:2.13} 
|\rho_{\ell}|\le 1 .
\end{equation} 

Now, we can estimate the remainder~\eqref{eq:2.5}. We obtain
\begin{equation}\label{eq:2.14} 
|\cR|\le R(\Delta^2) +2(\log\Delta)R(z\Delta^2) .
\end{equation} 

Next, we proceed to an estimation of our main term $W$. 
For our choice of $y_d$ given by~\eqref{eq:2.11} we find that
\begin{equation}\label{eq:2.15} 
y_p-y_{pd}=H^{-1}\min\bigl(\log p, \log\frac{\Delta}{d}\bigr)
\end{equation} 
and the formula~\eqref{eq:2.8}  becomes
\begin{equation}\label{eq:2.16} 
H^2W = K(\Delta) -\sum_{p<z}g(p)
\sum_{\substack{d<\Delta\\ (d,p)=1}}h(d)
\bigl\{\min \bigl(\log p, \log\frac{\Delta}{d}\bigr)\bigr\}^2
\end{equation} 
where, for any $u\le \Delta$
\begin{equation}\label{eq:2.17} 
K(u)= \sum_{d\le u} h(d)\bigl(\log\frac{\Delta}{d}\bigr)^2 . 
\end{equation} 
Note that, abbreviating $K=K(\Delta)$, we have $H^2\le J K$ 
by Cauchy's inequality. 

Having in mind that $g(p)$ is unstable at small primes we split the 
range of $p$ in~\eqref{eq:2.16} into two segments $p\le w$ and $w<p<z$ 
with some $2\le w<z$ at our disposal. Then we estimate $\min (.., ..)$ 
by $\log p$ and $\log(\Delta/d)$, respectively. We get
 \begin{equation}\label{eq:2.18} 
H^2W \ge K(\Delta) \bigl(1-\sum_{w<p<z}g(p)\bigr)-\sum_{p\le w}g(p)(\log p)^2
\sum_{\substack{d<\Delta\\ (d,p)=1}}h(d) .
\end{equation} 
We estimate the contribution of $p\le w$, $d\le w$, $(d,p)=1$ as follows 
(note that $g(p)h(d)\le h(pd)$). 
\begin{equation*}
\ssum_{\substack{p,d\le w \\ (p,d)=1}}h(pd)(\log p)^2
\le (\log w)\sum_{d\le w^2}h(d)\log d .
\end{equation*} 
If $\Delta\ge w^3$ we have 
$(\log w)(\log d)\le\frac12 (3\alpha\log \Delta/d)^2$, where
 \begin{equation}\label{eq:2.19} 
\alpha = (\log w)/\log \Delta .
\end{equation} 
Hence, this contribution is bounded by $\frac92 \alpha^2K(w^2)$ where, as 
in~\eqref{eq:2.17},  
 \begin{equation}\label{eq:2.20} 
K(w^2) = \sum_{d\le w^2}h(d)(\log \Delta/d)^2.
\end{equation} 
We drop the condition $(d,p)=1$ in the remaining sum 
 \begin{equation}\label{eq:2.21} 
J(w,\Delta)= \sum_{w<d\le\Delta}h(d).
\end{equation} 
Introducing the above estimates into~\eqref{eq:2.18}, we obtain
\begin{lem}
Let $z>w\ge 2$and $\Delta \ge w^3$. We have
 \begin{equation}\label{eq:2.22} 
H^2W\ge \bigl(1- \sum_{w< p \le z}g(p)\bigr)K 
-\tfrac92 \alpha^2K(w^2)-J(w, \Delta)G(w)
\end{equation} 
where 
 \begin{equation}\label{eq:2.23} 
G(w) = \sum_{p\le w}g(p)(\log p)^2.
\end{equation} 
\end{lem}

\newpage 

\section{\bf Two Assumptions}

So far, the inequality~\eqref{eq:2.22} holds without severe restrictions 
on the density functions $g(d)$, $h(d)$. To proceed further we accept two 
assumptions. 

{\bf Assumption 1.}    If $w$ is larger than some absolute constant, then
 \begin{equation}\label{eq:3.1} 
G(w)\le \tfrac32 (\log w)^2.
\end{equation} 

Note that if $g(p)p$ fluctuates in the interval $[0,2]$, then 
 \begin{equation*}
G(w)\le \sum_{p\le w}\tfrac{2}{p}(\log p)^2 =(\log w)^2 + O(\log w), 
\end{equation*}
so in this case we need $w$ to be sufficiently large to make the error 
term $O(\log w)$ strictly smaller than $\frac12 (\log w)^2$.

Let $K(w^2,\Delta)= K(\Delta) -K(w^2)$ denote the part of~\eqref{eq:2.17} 
that is complementary to~\eqref{eq:2.20};
  \begin{equation}\label{eq:3.2} 
K(w^2,\Delta)= \sum_{w^2<d\le \Delta}h(d)(\log \Delta/d)^2.
\end{equation} 
{\bf Assumption 2.} For $\Delta \ge w^3$ we have
  \begin{equation}\label{eq:3.3} 
J(w,\Delta) (\log \Delta)^2\le 3K(w^2,\Delta).
\end{equation} 

We shall illustrate how to verify~\eqref{eq:3.3} in special circumstances. 
But first we enjoy using both assumptions. By~\eqref{eq:3.1} 
and~\eqref{eq:3.3} we have 
  \begin{equation}\label{eq:3.4} 
J(w,\Delta) G(w)\le \tfrac92 \alpha^2 K(w^2,\Delta).
\end{equation} 
Recall that $\alpha=\log w/\log\Delta \le\tfrac13$. Hence Lemma 2.1 yields 
\begin{equation}\label{eq:3.5} 
H^2W\ge \nu K
\end{equation} 
where
\begin{equation}\label{eq:3.6} 
\nu = 1- \sum_{w<p\le z}g(p) -\tfrac92 \alpha^2 .
\end{equation} 
This lower bound is only interesting if $\nu$ is positive. Then $W$ is 
positive and we can apply the inequality $H^2\le J K$, getting 
\begin{theorem} Let $\Delta\ge w^3$ and $g(d)$ be such that~\eqref{eq:3.1}
and~\eqref{eq:3.3} hold. If $z>w$ and $\nu >0$, then
\begin{equation}\label{eq:3.7} 
J W\ge \nu .
\end{equation} 
\end{theorem}

{\bf Remarks} One is unlikely to obtain $\nu >0$ in normal situations. 
However, this can happen in exceptional circumstances, as for $g(p)$ given 
by~\eqref{eq:1.14} with an excepional real character $\chi (\mod q )$. These 
notes were designed mainly to handle such exceptional cases. 

\section{\bf Verification of Assumption 2}

Although $h(p)p$ is unpredictable at small primes, $h(d)d$ can be quite 
regular for large $d$ in the sense that
\begin{equation}\label{eq:4.1} 
\sum_{d\le x}h(d)d = cx +O\bigl( q^{\frac14} x^{\frac34}\bigr)
\end{equation} 
holds for every $x\ge 2$ with some constants $c>0$, $q\ge 2$ and an 
absolute implied constant in the error term. Having this formula, we can 
establish~\eqref{eq:3.3} by asymptotic evaluation of both sides. 

Recall that $\alpha =\log w/\log \Delta \le \tfrac13$. On the left side we get
\begin{equation*}
\begin{aligned}
J (w,\Delta)  & =\int_w^{\Delta}x^{-1}d\bigl(cx +O( q^{\frac14} x^{\frac34})\bigr)\\
& =c\log\tfrac{\Delta}{w} +O( q^{\frac14} w^{-\frac14})
=(1-\alpha)c\bigl(\log\Delta   +O(1)\bigr)
\end{aligned}
\end{equation*} 
provided that $w\gg qc^{-4}$. On the right side we get
\begin{equation*}
\begin{aligned}
K (w^2,\Delta)  & =\int_{w^2}^{\Delta}x^{-1}\bigl(\log \tfrac{\Delta}{x}\bigr)^2
d\bigl(cx +O( q^{\frac14} x^{\frac34})\bigr)\\
& =\tfrac{c}{3}\bigl( (\log\Delta)^3 -(\log w^2)^3\bigr) +
O\bigl( q^{\frac14} w^{-\frac14}(\log\Delta)^2\bigr)\\
& =\tfrac{c}{3}(1-8\alpha^3)(\log\Delta)^2\bigl(\log\Delta   +O(1)\bigr) .
\end{aligned}
\end{equation*} 
Hence, the ratio of the right side to the left side of ~\eqref{eq:3.3} is
\begin{equation*}
\tfrac{1-8\alpha^3}{1-\alpha}
+O\bigl(\tfrac{1}{\log \Delta}\bigr)\ge 1+\tfrac{\alpha}{6}
+O\bigl(\tfrac{1}{\log \Delta}\bigr) > 1 
\end{equation*}
provided that $w$ is sufficiently large to compensate for 
the error term $O(1/\log\Delta)$. This proves~\eqref{eq:3.3} if 
$\Delta\ge w^3$, $w\gg qc^{-4}$. 

\smallskip

{\bf Example} Let $g(p)p$ be given by~\eqref{eq:1.14}. Then
\begin{equation}\label{eq:4.2} 
h(p)p = \bigl(1-\tfrac1p\bigr)^{-1}\bigl(1-\tfrac{\chi(p)}{p}\bigr)^{-1}
\bigl(1 +\chi(p)(1-\tfrac{1}{p})\bigr) 
\end{equation}
if $p\nmid q$ and $h(p)=0$ if $p|q$. 
Hence, the series
\begin{equation*}
D(s)=\sum_dh(d)d^{1-s} =\prod_p\bigl(1+h(p)p^{1-s}\bigr)
=\zeta(s)L(s,\chi)E(s)
\end{equation*}
has analytic continuation to ${\rm Re} s >\tfrac12$. The local factors 
\begin{equation*}
\begin{aligned}
E_p(s) & = \bigl(1-p^{-s}\bigr)\bigl(1-\chi(p)p^{-s}\bigr)
\bigl(1+h(p)p^{1-s}\bigr)\\
& = 1+a_1p^{-s-1} +a_2p^{-2s} +a_3p^{-3s} 
\end{aligned}
\end{equation*} 
have $a_1$, $a_2$, $a_3$ bounded. Hence~\eqref{eq:4.1} follows by standard 
contour integration with the constant
\begin{equation*}
c =\res_{s=1}D(s) =L(1,\chi)E(1) .
\end{equation*} 
We have 
\begin{equation*}
 1+h(p)=\bigl(1-g(p)\bigr)^{-1}=\bigl(1-\tfrac1p\bigr)^{-1}
\bigl(1-\tfrac{\chi(p)}{p}\bigr)^{-1} 
\end{equation*} 
if $p\nmid q$ and $1$ if $p|q$. Hence $E(1)=\varphi(q)/q$ so the residue is 
\begin{equation}\label{eq:4.3} 
c=L(1,\chi)\varphi(q)/q\gg \varphi(q)q^{-\tfrac32}(\log q)^{-2},
\end{equation} 
by the Dirichlet class number formula.  Therefore, we have proved 
that~\eqref{eq:3.3} holds as long as
\begin{equation}\label{eq:4.4} 
\Delta \ge w^3, \quad w\ge q^3, 
\end{equation} 
and $q$ is sufficiently large.

\section{\bf Exceptional Primes}

We are going to apply the theorem to estimate the sum of $a_p$ over primes 
$p\le x$. In the context of these notes we say that the density function 
$g(d)$ is ``exceptional'' if 
\begin{equation}\label{eq:5.1} 
\delta(w,z)=  \sum_{ w<p\le z}g(p)
\end{equation} 
is bounded, smaller than $1$. Specifically, we take $\Delta=w^3$, 
$z=\sqrt x\ge w^3$ and we make the following

{\bf Assumption 3.}  
\begin{equation}\label{eq:5.2} 
\delta = (w,\sqrt x) \le \tfrac14 .
\end{equation} 

\smallskip
We then have
\begin{equation}\label{eq:5.3} 
\nu = 1-\delta (w,\sqrt x) -\tfrac12 \ge \tfrac14 .
\end{equation} 
so that~\eqref{eq:3.7} implies $W\ge \tfrac14 V$, where
\begin{equation}\label{eq:5.4} 
V=\prod_{p<w^3}\bigl(1-g(p)\bigr) . 
\end{equation} 
If $\cA=(a_n)$ is supported on $\sqrt x <n \le x$, then
\begin{equation}\label{eq:5.5} 
\sum_pa_p = S(\cA,z)\ge S^-(\cA,z) =XW +\cR  
\end{equation} 
where the remainder $\cR$ is bounded by (see~\eqref{eq:1.5} and~\eqref{eq:2.14})
\begin{equation}\label{eq:5.6} 
\cR \ll R \log x = (\log x)\sumb_{d<\sqrt x w^6}\tau_3(d)|r_d| .
\end{equation} 
Hence we conclude the following
\begin{prop}
If $g$ is exceptional, that is~\eqref{eq:5.2} holds, then
\begin{equation}\label{eq:5.7} 
\sum_{\sqrt x <p\le x}a_p \ge \tfrac14 XV - R\log x .    
\end{equation} 
\end{prop}
Recall that we work under Assumptions 1, 2, 3. The first two assumptions 
are verified in the case of $g(p)$ given by~\eqref{eq:1.14}. The required 
conditions are those in~\eqref{eq:4.4}. We choose
\begin{equation}\label{eq:5.8} 
w=q^3 .
\end{equation} 
We obtain
\begin{equation*}
V=\prod_{\substack{p<w^3\\p\nmid q}}\bigl(1-\frac1p \bigr)
\bigl(1- \frac{\chi(p)}{p}\bigr)>\tfrac{1}{21} V(q)\frac{q}{\varphi(q)}
\end{equation*} 
where
\begin{equation}\label{eq:5.9} 
V(q)= \prod_{p<q^2}\bigl(1-\frac1p \bigr)
\bigl(1- \frac{\chi(p)}{p}\bigr) .
\end{equation} 

\section{\bf The Exceptional Character}

Let $\lambda = 1 * \chi$, so $\lambda (p) = 1+\chi(p)$. We 
have proved in (24.20) of [Opera] that 
\begin{equation}\label{eq:6.1} 
\sum_{q^3 < p\le \sqrt x}\lambda(p)p^{-1}< (1-\beta)\log x +O(q^{-3/4})
\end{equation} 
if $x\ge q^6$, where $\beta$ is any real zero of $L(s,\chi)$ and the 
implied constant is absolute. Note that 
\begin{equation}\label{eq:6.2} 
1-\beta \gg q^{-1/2} (\log q)^{-2},
\end{equation} 
so the error term in the inequality~\eqref{eq:6.1} is negligible. 
Now, our exceptional condition~\eqref{eq:5.2} holds in the segment
\begin{equation}\label{eq:6.3} 
q^6\le x\le e^{1/4(1-\beta)}
\end{equation} 
which is non-empty if
\begin{equation}\label{eq:6.4} 
(1-\beta)\log q\le \tfrac{1}{24}. 
\end{equation} 
Therefore, we shall say that the real primitive character $\chi$ is 
exceptional if $L(s,\chi)$ has a real zero $\beta$ satisfying ~\eqref{eq:6.4}. 
If also $\chi(a) = 1$ its effect pulls in an unfavourable direction.  

Under this exceptional situation we are able to estimate the least prime 
in an arithmetic progression,
\begin{equation}\label{eq:6.5} 
p\equiv a (\mod q) \quad {\rm with}\,\,\, \chi(a)=1 .
\end{equation} 
To this end, consider the sequence
\begin{equation}\label{eq:6.6} 
a_n= \lambda(n)
\end{equation} 
for $ n \equiv a (\mod q)$, $1\le n \le x$. We set $a_n=0$ otherwise. We have 
shown in Section 24.2 of [Opera] that $\cA$ has our sieve properties 
with the density function~\eqref{eq:1.14} and 
\begin{equation}\label{eq:6.7} 
X=2L(1,\chi)xq^{-1} .
\end{equation} 
The individual error terms satisfy 
\begin{equation}\label{eq:6.8} 
r_d \ll \tau_3(d)\sqrt{x/d}
\end{equation} 
if $d\le x$, where the implied constant is absolute. Hence, our 
remainder~\eqref{eq:5.6} satisfies 
\begin{equation}\label{eq:6.9} 
R \log x \ll w^3x^{\tfrac 34}(\log x)^{10} .
\end{equation} 
We want this bound to be insignificant by comparison with the main term
\begin{equation}\label{eq:6.10} 
\tfrac14 XV \gg xq^{-\tfrac32}(\log q)^{-2} .
\end{equation} 
This is the case if $x\ge q^{43}$. We obtain

\begin{cor}Suppose $L(s,\chi)$ has a real zero $\beta$ with
\begin{equation}\label{eq:6.11} 
(1-\beta)\log q \le \tfrac{1}{172}. 
\end{equation} 
Let $\chi(a)=1$. Then, for $q^{43}\le x \le e^{1/4(1-\beta)}$, we have 
\begin{equation}\label{eq:6.12} 
\pi(x;q,a)\ge L(1,\chi)\frac{V(q)}{\varphi(q)}\frac{x}{168} ,
\end{equation} 
where $V(q)$ is the product~\eqref{eq:5.9}. 
\end{cor}

\end{document}